\documentclass[twocolumn,amsmath,amssymb,aps,secnumarabic,%
    nofootinbib]{revtex4-1}
\usepackage{times,mathptmx,amsthm,tikz}
\usetikzlibrary{decorations.markings}

\newtheorem{theorem}{Theorem}

\newcommand{\R}{\mathbb{R}}
\newcommand{\C}{\mathbb{C}}

\newcommand{\fig}[1]{Figure~\ref{#1}}
\newcommand{\eq}[2]{\begin{equation}\label{#1}#2\end{equation}}

\newenvironment{fullfigure}[2]
    {\begin{figure}[htb]\begin{center}\def\fullfiga{#1}\def\fullfigb{#2}}
    {\vspace{\baselineskip}\caption{\fullfigb.}\label{\fullfiga}
    \end{center}\end{figure}}

\begin{document}
\title{Norms as a function of $p$ are linearly independent
    in finite dimensions}

\author{Greg Kuperberg}
\affiliation{University of California, Davis}

\begin{abstract} We show that there are no non-trivial linear dependencies
among $p$-norms of vectors in finite dimensions that hold for all $p$.
The proof is by complex analytic continuation.
\end{abstract}
\maketitle

\begin{theorem} Let $v_1,v_2,\ldots,v_n$ be non-zero vectors with
$v_k \in \R^{d_k}$.  Suppose that
\eq{e:main}{\alpha_1||v_1||_p + \alpha_2||v_2||_p +
    \cdots + \alpha_n||v_n||_p = 0}
for all $p \in [a,b]$ with $1 \le a < b \le \infty$.  Then the equation is
trivial in the following sense:  Call two of the vectors equivalent if they
differ by adding zeros, permuting or negating coordinates, and rescaling.
Then the terms of \eqref{e:main} in each equivalence class, with the given
coefficients, already sum to zero.
\end{theorem}

The result was stated as a question by Steve Flammia on MathOverflow.

\begin{proof} Suppose by induction that $n$ has the least value.  We can
write
\eq{e:expln}{||v_k||_p = \exp \left(\frac{\ln \left[\exp(\beta_1 p) +
    \exp(\beta_2 p) + \cdots + \exp(\beta_d p)\right]}{p}\right),}
by letting $\beta_j = \ln |v_{k,j}|$, ignoring vanishing coordinates,
and taking $d = d_k$.  We will consider the complex analytic continuation
of $||v_k||_p$ and equation \eqref{e:main} using this formula.  (The
technique of complex analytic continuation in $p$ was also used by Thorin
\cite{Thorin:convexity} for a different purpose.)

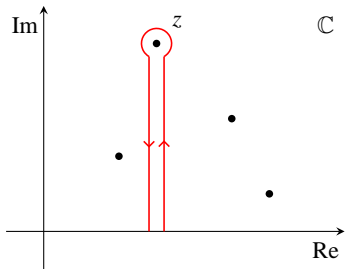
\begin{fullfigure}{f:continuation}{Analytic continuation of $||v_k||_p$
    around a zero $z$ of $f_k(p)$}
\tikzset{midto/.style={postaction={decorate,
    decoration={markings,mark=at position .5 with
    {\draw (-.035,-.07) -- (.035,0) -- (-.035,.07);}}}}}
\colorlet{darkred}{red!75!black}
\begin{tikzpicture}[>=stealth]
\draw[->] (-.5,0) -- (4,0);
\draw[->] (0,-.5) -- (0,3);
\fill (1,1) circle (.05);
\fill (1.5,2.5) circle (.05);
\fill (2.5,1.5) circle (.05);
\fill (3,.5) circle (.05);
\draw[red,midto,semithick] (1.6,0) -- (1.6,2.3268);
\draw[red,semithick] (1.6,2.3268) arc (-60:240:.2);
\draw[red,midto,semithick] (1.4,2.3268) -- (1.4,0);

\draw (1.6,2.6) node[anchor=south west] {$z$};
\draw (-.25,2.75) node {$\mathrm{Im}$};
\draw (3.75,-.25) node {$\mathrm{Re}$};
\draw (3.75,2.75) node {$\C$};
\end{tikzpicture}
\end{fullfigure}

Let $f_k(p)$ be the exponential sum in the argument of the logarithm in
equation \eqref{e:expln}. It is also called an exponential polynomial,
because it is a polynomial in $\exp(p)$ with non-integer exponents
$\beta_1,\beta_2,\ldots,\beta_d$.  Plainly $f_k(p)$ is an entire function
of $p$, which implies that it has isolated zeros.  Also, none of its zeros
are on the real axis.  Thus $||v_k||_p$ is multivalued analytic
on all of $\C$, except for $p=0$ and except for the zeros of $f_k(p)$.

As in \fig{f:continuation}, we can continue $f_k(p)$ along a path that
begins and ends on the positive real axis and travels around a small loop
around an $m$-fold zero $z$.  This continuation multiplies $||v_k||_p$
by $\exp(2\pi im/z)$, which is distinct for all $m$ because we must
have $\mathrm{Im}\ z \ne 0$.  This changes equation \eqref{e:main} to create a
new linear dependence of the same form.  This contradicts the minimality
of $n$, unless the new equation is proportional to the old one.  Thus,
each $f_k(p)$ has the same multiset of zeroes.

By the Hadamard factorization theorem \cite[\S 3.4]{Conway:gtm}, if $f(p)$
is an entire function with $|f(p)| = \exp(O(|p|))$, then it is determined
by its zeroes up to a factor of $a \exp(\beta p)$.  (This step is due to
the referee.  The author's first proof used a structure theorem for
exponential polynomials due to Ritt \cite{Ritt:exp}.)  Thus for any two
$j \ne k$, we can write
\eq{e:ratio}{f_j(p) = a\exp(\beta p)f_k(p),}
or
$$||v_j||_p = a^{1/p}\exp(\beta)||v_k||_p.$$
This also implies more than one linear dependence of the form \eqref{e:main}
unless $a=1$ for all pairs of vectors.  Finally, equation \eqref{e:ratio}
with $a=1$ is only possible if the vectors $v_j$ and $v_k$ are equivalent,
given that exponential monomials with distinct exponents are linearly
independent.
\end{proof}

We leave open the question of finding solutions to equation \eqref{e:main}
in infinite-dimensional $\ell^p$ or $L^p$ spaces.

\acknowledgments

The author would like to thank Steve Flammia for raising the question,
and him, Mark Meckes, and the referee for useful comments.  The author
was partly supported by NSF grants DMS-0606795 and CCF-1013079.


\providecommand{\bysame}{\leavevmode\hbox to3em{\hrulefill}\thinspace}
\providecommand{\MR}{\relax\ifhmode\unskip\space\fi MR }
\providecommand{\MRhref}[2]{%
  \href{http://www.ams.org/mathscinet-getitem?mr=#1}{#2}
}
\providecommand{\href}[2]{#2}
\providecommand{\eprint}{\begingroup \urlstyle{tt}\Url}

\end{document}